\UseRawInputEncoding
\documentclass[reqno,10.5pt]{amsart}

\usepackage{amsthm,amsmath,amssymb}
\usepackage{mathrsfs,amsfonts,dsfont,functan,extarrows,mathtools}
\usepackage[colorlinks]{hyperref}

\usepackage{marginnote}
\usepackage{xcolor}

\makeatletter

\newcommand{\Rmnum}[1]{\expandafter\@slowromancap\romannumeral #1@}
\makeatother

\def\XXint#1#2#3{{\setbox0=\hbox{$#1{#2#3}{\int}$}	
\vcenter{\hbox{$#2#3$}}\kern-.5\wd0}}

\newtheorem{theorem}{Theorem}[section]
\newtheorem{proposition}[theorem]{Proposition}

\newtheorem{lemma}[theorem]{Lemma}

\numberwithin{equation}{section}
\allowdisplaybreaks

\newcounter{wronumber}

\begin{document}

\title[]{Gradient estimates for Donaldson's equation on a compact K\"ahler manifold}

\author[L. Zhang]{Liangdi Zhang}
\address[Liangdi Zhang]{\newline 1. Yanqi Lake Beijing Institute of Mathematical Sciences and Applications, Beijing 101408, P. R. China
\newline 2. Yau Mathematical Sciences Center, Tsinghua University, Beijing 100084, P. R. China}
\email{ldzhang91@163.com}



\begin{abstract}
We prove a gradient estimate for Donaldson's equation
\[\omega\wedge(\chi+\sqrt{-1}\partial\overline{\partial}\varphi)^{n-1}=e^F(\chi+\sqrt{-1}\partial\overline{\partial}\varphi)^n\]
(and its parabolic analog) on an $n$-dimensional compact K\"ahler manifold $(M,\omega)$ with another Hermitian metric $\chi$ directly from the uniform upper bounds for $tr_\omega\chi_\varphi$ and Alexandrov-Bakelman-Pucci (ABP) maximum principle.

\vspace*{5pt}


\end{abstract}

\maketitle


\section{Introduction} 
\label{sec:1}
Let $(M,\omega)$ be a compact K\"ahler manifold of complex dimension $n$ and $\chi$ be another K\"ahler metric on $M$. In local coordinates,
\[\omega=\sqrt{-1}g_{i\bar{j}}dz^i\wedge d\bar{z}^j\ \text{and}\ \chi=\sqrt{-1}\chi_{i\bar{j}}dz^i\wedge d\bar{z}^j.\]
$\mathcal{H}_\chi$ denotes the set of all $\varphi\in C^\infty(M,\mathbb{R})$ satisfying
\[\chi_\varphi:=\chi+\sqrt{-1}\partial\overline{\partial}\varphi>0.\]

The $J$-flow introduced by S. K. Donaldson \cite{do99} and X. Chen \cite{chen00} independently is a parabolic flow of
\begin{equation}\label{1.4}
\frac{\partial\varphi}{\partial t}=c-\frac{\omega\wedge\chi_\varphi^{n-1}}{\chi_\varphi^n}
\end{equation}
with $c$ is the constant given by
\[c=\frac{\int_M\omega\wedge\chi^{n-1}}{\int_M\chi^n}.\]


X. Chen \cite{sw7} proved the $J$-flow \eqref{1.4} exists for all time with any smooth initial data. It is known that the critical point to \eqref{1.4} is given by the Donaldson's equation
\begin{equation}\label{1.5}
\omega\wedge\chi^{n-1}=c\chi^n,
\end{equation}
where the K\"ahler metric $\chi$ is in the cohomology class $[\chi]$ (see \cite{do99,chen00} for details). S. K. Donaldson \cite{do99}  found that if \eqref{1.5} has a solution in $[\chi]$ then
\begin{equation}\label{1.6}
[nc\chi]-[\omega]>0,
\end{equation}
and remarked a natural conjecture that the necessary condition \eqref{1.6} is also sufficient for the existence for a solution to \eqref{1.5} in $[\chi]$.
Donaldson's conjecture was proved by X. Chen \cite{sw7} in an elliptic way while by B. Weinkove \cite{wein04} in a parabolic way for the case of $n=2$, while by B. Weinkove \cite{wein06} in higher dimensions. Moreover, J. Song and B. Weinkove \cite{sw08} showed that \eqref{1.5} has a smooth solution in $[\chi]$ if and only if there exists a metric $\chi'\in[\chi]$ with
\begin{equation}\label{1.7}
(nc\chi'-(n-1)\omega)\wedge\chi'^{n-2}>0,
\end{equation}
and \eqref{1.7} is also a sufficient and necessary condition for the $J$-flow \eqref{1.4} converges in $C^\infty$ to a function $\phi_\infty\in\mathcal{H}_\chi$ with $\chi_{\phi_\infty}$ satisfies \eqref{1.5}. In 2021, G. Chen \cite{gchen2021} considered the $J$-equation
\begin{equation}\label{1.17}
tr_{\omega_\varphi}\chi=c
\end{equation}
and proved that there exists a K\"ahler metric $\omega_\varphi=\omega+\sqrt{-1}\partial\bar{\partial}\varphi>0$ satisfying the $J$-equation \eqref{1.17} if and only if $(M,[\omega],[\chi])$ is uniformly $J$-stable.

On a compact Hermitian manifold, V. Tosatti and B. Weinkove \cite{tw10} obtained a uniform $C^\infty$ estimate for a smooth solution $\varphi$ to the complex Monge-Amp\`{e}re equation
\begin{equation}\label{1.0}
(\omega+\sqrt{-1}\partial\overline{\partial}\varphi)^n=e^F\omega^n
\end{equation}
with $\omega+\sqrt{-1}\partial\overline{\partial}\varphi>0$ and proved the Hermitian version of Calabi conjecture that every representative of the first Bott-Chern class can be
represented as the first Chern form of a Hermitian metric of the form $\omega+\sqrt{-1}\partial\overline{\partial}\varphi$. In 2011, M. Gill \cite{gill11} proved long-time existence to the parabolic complex Monge-Amp\`{e}re equation
\begin{equation}\label{1.8}
\frac{\partial\varphi}{\partial t}=\log\frac{(\omega+\sqrt{-1}\partial\overline{\partial}\varphi)^n}{\omega^n}-F
\end{equation}
with $\omega+\sqrt{-1}\partial\overline{\partial}\varphi>0$ and the initial condition $\varphi(\cdot,0)=0$ on a compact Hermitian manifold $(M,\omega)$ of complex dimension $n$, and provided a parabolic proof of the main theorem in
V. Tosatti and B. Weinkove \cite{tw10}. More generally, W. Sun \cite{sun15} derived uniform $C^\infty$ estimates for the normalized solution and the $C^\infty$ convergence of the parabolic complex Monge-Amp\`{e}re type equations
\begin{equation}\label{1.9}
\frac{\partial u}{\partial t}=\log\frac{\chi_u^n}{\chi_u^{n-\alpha}\wedge\omega^\alpha}-\log\psi\ (1\leq \alpha\leq n)
\end{equation}
with $\chi$ is a smooth positive $(1,1)$ form on a compact Hermitian manifold $(M^n,\omega)$.

The second order estimate for the complex Monge-Amp\`{e}re equation \eqref{1.0} on a compact K\"ahler manifold plays a connecting link role in S. T. Yau's original proof of the Calabi conjecture (see \cite{yau} for details). It is noted that the second order estimate in S. T. Yau \cite{yau} is directly obtained once the zeroth order estimate is known and does not make use of gradient estimates. However, there has been an increasing interest for researchers to study gradient estimates on the complex Monge-Amp\`{e}re equation on compact manifolds. Z. B{\l}ocki \cite{blocki09} derived a gradient estimate for the complex Monge-Amp\`{e}re equation on a compact K\"ahler manifold directly from its zeroth order estimate, without using the second order estimates. B. Guan and Q. Li \cite{GL09} generalized this estimate to a complex Monge-Amp\`{e}re equation satisfying a Dirichlet condition on a compact Hermitian manifold with smooth boundary. Moreover, X. Chen and W. He \cite{ch12} improved the gradient estimate by Z. B{\l}ocki \cite{blocki09} that the right-hand-side of \eqref{1.0} has rather weaker regularity. More recently, X. Chen and J. Cheng \cite{cc19} showed $L^p$ ($p<\infty$) and $L^\infty$ estimates for the gradient of the solutions to the complex Monge-Amp\`{e}re equation \eqref{1.0} in terms of the continuity of $F$ and without assuming any bounds on its derivatives. In 2021, B. Guo, D. H. Phong and F. Tong \cite{gpt21} proved a gradient estimate for the complex Monge-Amp\`{e}re equation \eqref{1.0} on a compact K\"ahler manifold with boundary by using the Alexandrov-Bakelman-Pucci (ABP) maximum principle.

X. Zhang \cite{xzhang10} obtained gradient estimates for geometric solutions of the complex Hessian equations
\begin{equation}\label{1.11}
\sigma_k(g^{-1}(g_{i\bar{j}}+\phi_{i\bar{j}}))=f
\end{equation}
and
\begin{equation}\label{1.12}
\sigma_k(g_{i\bar{j}}+\phi_{i\bar{j}}+\mu\phi_i\phi_{\bar{j}})=f
\end{equation}
on a compact Hermitian manifold $(M,g_{i\bar{j}})$ by maximum principle arguments, where $\sigma_k$ is the $k$-th $(1\leq k\leq n)$ elementary symmetric function and $\mu$, $f$ $(f>0)$ are smooth functions on $M$. It is clear that \eqref{1.11} and \eqref{1.12} are natural generalizations of the complex Monge-Amp\`{e}re equation. In 2018, R. Yuan \cite{ryuan18} derived gradient estimates from the Bernstein
method and other ordered apriori estimates on a class of second order fully nonlinear elliptic equations containing gradient terms on compact Hermitian manifolds under some assumptions. See also \cite{twwy23,sjdg20,sjdg21,sjdg23,sjdg25,sjdg35,sjdg36,sjdg38,sjdg,twwy64,ryuan22} for related work on apriori estimates for fully nonlinear equations on manifolds.

In this paper, we consider gradient estimates on smooth solutions  $\varphi\in\mathcal{H}_\chi$ to the Donaldson's equation
\begin{equation}\label{1.1}
\omega\wedge\chi_\varphi^{n-1}=e^F\cdot\chi_\varphi^n
\end{equation}
and the parabolic Donladson's equation
\begin{equation}\label{1.10}
\frac{\partial \varphi}{\partial t}=\log\frac{\chi_\varphi^n}{\omega\wedge\chi_\varphi^{n-1}}+F,
\end{equation}
respectively, on an $n$-dimensional compact K\"ahler manifold $(M,\omega)$ with $\chi$ is another Hermitian metric on $M$ and $F\in C^\infty(M,\mathbb{R})$ is a given smooth function. It is noted that the Donaldson's equation \eqref{1.1} is a generalization of the complex Monge-Amp\`{e}re equation \eqref{1.0} while the parabolic Donladson's equation \eqref{1.10} is a special case of \eqref{1.9}. The main analytic tool we apply is ABP maximum principle. It is noted that the usage of ABP maximum principle for the complex Monge-Amp\`{e}re equation was originated in work of S. Y. Cheng and S. T. Yau. More recently, ABP maximum principle was revisited by Tosatti-Weinkove \cite{tw21}, Guo-Phong-Tong \cite{gpt217,gpt21} and Chen-Cheng \cite{cc21} in the study of complex Monge-Amp\`{e}re equations, by Chen-Cheng \cite{chencheng1} in cscK equations and by G. Sz\'ekelyhidi \cite{sjdg} in fully nonlinear elliptic equations. We do not use any known apriori estimates other than second order estimates.

In 2014, Y. Li \cite{yli14} proved $C^\infty$ estimate for the Donaldson's equation \eqref{1.1} when $\omega$ is a Hermitian metric and $\chi-\frac{n-1}{ne^F}\omega>0$. The first goal of this paper is to complete apriori estimates of \eqref{1.1} in the case of $\omega$ is K\"ahler. To be precise, we establish a gradient estimate from the second order estimate by Y. Li \cite{yli14}.
\begin{theorem}\label{thm1.1}
Let $(M,\omega)$ be an $n$-dimensional compact K\"ahler manifold (without boundary) and $\chi$ be another Hermitian metric on $M$. Assume that
\begin{equation}\label{1.2}
\chi-\frac{n-1}{ne^F}\omega>0.
\end{equation}
If $\varphi\in\mathcal{H}_\chi$ is a solution to the Donaldson's equation \eqref{1.1} on $M$, then there exists a uniform constant $C$ depending on $M$, $\omega$, $\chi$ and $F$, and a positive constant $\lambda$ depending only on $\omega$ and $\chi$ such that
\begin{equation}\label{1.3}
|\nabla\varphi|_\omega^2\leq Ce^{\lambda(\varphi-\inf_M\varphi)}
\end{equation}
on $M$.
\end{theorem}

Following the definition given in \cite{sun15}, we say the Hermitian metric $\chi$ satisfies the cone condition if $[\chi]\in\mathcal{C}(F)$, where
\[\mathcal{C}(F):=\{[\chi]:\ \text{there exists}\ \chi'\in[\chi]^+,\ n\chi'^{n-1}>(n-1)e^{-F}\chi'^{n-2}\wedge\omega\}\]
with $[\chi]^+:=\{\chi'\in[\chi]:\chi'>0\}$. The second goal in this paper is to prove a gradient estimate for the parabolic Donaldson's equation \eqref{1.10} by only using the second order estimate in W. Sun \cite{sun15} and the ABP maximum principle of parabolic version.
\begin{theorem}\label{thm1.2}
Let $(M,\omega)$ be an $n$-dimensional compact K\"ahler manifold (without boundary) and $\chi$ be another Hermitian metric on $M$. Assume that $[\chi]\in\mathcal{C}(F)$ and
\[e^F\chi^n\leq\omega\wedge\chi^{n-1}.\]
If $\varphi$ is a solution to the parabolic Donaldson's equation \eqref{1.10} on $M\times[0,T]$ for $0<T<\infty$, then there exists a uniform constant $\bar{C}$ depending only on the initial geometric data and $T$, and a positive constant $\lambda$ depending only on $\omega$ and $\chi$ such that
\begin{equation}\label{1.3}
|\nabla\varphi|_\omega\leq \bar{C} e^{\lambda(\varphi-\inf_{M\times[0,T]}\varphi)}
\end{equation}
on $M\times[0,T]$.
\end{theorem}

In the rest of this paper, we prove Theorem \ref{thm1.1} and Theorem \ref{thm1.2} in Section \ref{sec:2} and Section \ref{sec:3}, respectively.
\section{Gradient estimate for Donaldson's equation} 
\label{sec:2}
Motivated by the work of Guo-Phong-Tong \cite{gpt21}, we estimate upper bounds for $|\nabla\varphi|_\omega$ with $\varphi\in\mathcal{H}_\varphi$ solving \eqref{1.1}. As in \cite{yli14}, we define a Hermitian metric $h$ on $M$ with its inverse matrix is locally given by
\[h^{i\bar{j}}=\chi_\varphi^{i\bar{l}}\chi_\varphi^{k\bar{j}}g_{k\bar{l}}.\]

Let $R(\omega)_{i\bar{j}k\bar{l}}$ be the curvature tensor of $\omega$. $\Delta_\omega$ and $\Delta_h$ denote the Laplacian operator associated to the K\"ahler metric $\omega$ and the Hermitian metric $h$, respectively.

First of all, we calculate a formula for $\Delta_h|\nabla\varphi|_\omega^2$.

\begin{proposition}\label{prp2.1}
The following formula holds on $M$.
\begin{eqnarray}\label{2.1}
\Delta_h|\nabla\varphi|_\omega^2&=&-2ne^F\nabla F\cdot_\omega \nabla\varphi-2tr_h\nabla\chi\cdot_\omega\nabla\varphi\notag\\
&&+2h^{i\bar{j}}g^{k\bar{q}}g^{p\bar{l}}R(\omega)_{i\bar{j}k\bar{l}}\varphi_p\varphi_{\bar{q}}\notag\\
&&+h^{i\bar{j}}g^{k\bar{l}}(\varphi_{ki}\varphi_{\bar{j}\bar{l}}+\varphi_{k\bar{j}}\varphi_{i\bar{l}}).
\end{eqnarray}
\end{proposition}

\begin{proof}
Without loss of generality, we calculate in a holomorphic normal coordinate system for the K\"ahler metric $\omega$ at an arbitrary point $x$. By direct computation, we have
\begin{eqnarray}\label{2.2}
\Delta_h|\nabla\varphi|_\omega^2&=&h^{i\bar{j}}(g^{k\bar{l}}\varphi_k\varphi_{\bar{l}})_{\bar{j}i}\notag\\
&=&h^{i\bar{j}}(g^{k\bar{l}},_{\bar{j}}\varphi_k\varphi_{\bar{l}}+g^{k\bar{l}}\varphi_{k\bar{j}}\varphi_{\bar{l}}+g^{k\bar{l}}\varphi_k\varphi_{\bar{l}\bar{j}})_i\notag\\
&=&h^{i\bar{j}}(g^{k\bar{l}},_{\bar{j}i}\varphi_k\varphi_{\bar{l}}+\varphi_{k\bar{j}i}\varphi_{\bar{k}}+\varphi_{k\bar{j}}\varphi_{\bar{k}i}\notag\\
&&+\varphi_{ki}\varphi_{\bar{k}\bar{j}}+\varphi_k\varphi_{\bar{k}\bar{j}i})\notag\\
\end{eqnarray}

At such a point, it is clear that
\begin{equation}\label{2.3}
g^{k\bar{l}},_{\bar{j}i}=R(\omega)_{i\bar{j}k\bar{l}},
\end{equation}
and
\begin{equation}\label{2.4}
\varphi_{\bar{k}\bar{j}i}=\varphi_{\bar{k}i\bar{j}}+R(\omega)_{i\bar{j}l\bar{k}}\varphi_{\bar{l}}=\varphi_{i\bar{j}\bar{k}}+R(\omega)_{i\bar{j}l\bar{k}}\varphi_{\bar{l}}.
\end{equation}

Applying \eqref{2.3} and \eqref{2.4} to \eqref{2.2} yields
\begin{eqnarray}\label{2.5}
\Delta_h|\nabla\varphi|_\omega^2&=&2h^{i\bar{j}}R(\omega)_{i\bar{j}k\bar{l}}\varphi_{\bar{k}}\varphi_l+h^{i\bar{j}}(\varphi_{i\bar{j}k}\varphi_{\bar{k}}+\varphi_k\varphi_{i\bar{j}\bar{k}})\notag\\
&&+h^{i\bar{j}}(\varphi_{k\bar{j}}\varphi_{\bar{k}i}+\varphi_{ki}\varphi_{\bar{k}\bar{j}}).
\end{eqnarray}

From \eqref{1.1}, we know that
\begin{equation}\label{2.5x}
ne^F=n\frac{\omega\wedge\chi_\varphi^{n-1}}{\chi_\varphi^n}=tr_{\chi_\varphi}\omega:=\chi_\varphi^{i\bar{j}}g_{i\bar{j}}.
\end{equation}

Hence,
\[ne^F\cdot F_k=-h^{i\bar{j}}\varphi_{i\bar{j}k}-h^{i\bar{j}}\chi_{i\bar{j},k}\ \text{and}\ ne^F\cdot F_{\bar{k}}=-h^{i\bar{j}}\varphi_{i\bar{j}\bar{k}}-h^{i\bar{j}}\chi_{i\bar{j},\bar{k}}.\]
Moreover, we obtain
\begin{eqnarray}\label{2.6}
h^{i\bar{j}}(\varphi_{i\bar{j}k}\varphi_{\bar{k}}+\varphi_k\varphi_{i\bar{j}\bar{k}})&=&-ne^F(F_k\varphi_{\bar{k}}+F_{\bar{k}}\varphi_k)-h^{i\bar{j}}(\chi_{i\bar{j},k}\varphi_{\bar{k}}+\chi_{i\bar{j},\bar{k}}\varphi_k)\notag\\
&=&-2ne^F\nabla F\cdot_\omega\nabla\varphi-2tr_h\nabla\chi\cdot_\omega\nabla\varphi.
\end{eqnarray}

\eqref{2.1} follows by applying \eqref{2.6} to \eqref{2.5} since it is independent of the choice of holomorphic normal coordinate systems.
\end{proof}

Let $-K$ be a uniform lower bound of the bisectional curvature of $\omega$ and $\lambda$ be a positive constant (depending only on $\omega$ and $\chi$) to be determined later.

\begin{lemma}\label{lem2.1}
Define $G:=e^{-\lambda(\varphi-\inf_M\varphi)}|\nabla\varphi|_\omega^2$, then there exists a uniform constant $C$ that depends only on $M$, $\omega$, $\chi$ and $F$ so that
\begin{equation}\label{2.7}
\Delta_hG\geq-CG^\frac{1}{2}+(\lambda tr_h\chi-2Ktr_h\omega)G-3\lambda ne^FG.
\end{equation}
\end{lemma}

\begin{proof}
Note that
\begin{equation}\label{2.9}
h^{i\bar{j}}(\chi_\varphi)_{i\bar{j}}=\chi_\varphi^{i\bar{l}}\chi_\varphi^{k\bar{j}}g_{k\bar{l}}(\chi_\varphi)_{i\bar{j}}=\chi_\varphi^{k\bar{l}} g_{k\bar{l}}=ne^F,\notag
\end{equation}
where we used \eqref{2.5x}. Then we have
\begin{eqnarray}\label{2.10}
-h^{i\bar{j}}\varphi_{i\bar{j}}&=&-h^{i\bar{j}}(\chi_\varphi)_{i\bar{j}}+h^{i\bar{j}}\chi_{i\bar{j}}\notag\\
&=&-ne^F+tr_h\chi.
\end{eqnarray}
Moreover,
\begin{eqnarray}\label{2.8}
\Delta_he^{-\lambda(\varphi-\inf_M\varphi)}&=&h^{i\bar{j}}(e^{-\lambda(\varphi-\inf_M\varphi)})_{\bar{j}i}\notag\\
&=&h^{i\bar{j}}(-\lambda\varphi_{\bar{j}}e^{-\lambda(\varphi-\inf_M\varphi)})_{i}\notag\\
&=&h^{i\bar{j}}(-\lambda\varphi_{i\bar{j}}+\lambda^2\varphi_i\varphi_{\bar{j}})e^{-\lambda(\varphi-\inf_M\varphi)}\notag\\
&=&(-\lambda ne^F+\lambda tr_h\chi+\lambda^2|\nabla\varphi|_h^2)e^{-\lambda(\varphi-\inf_M\varphi)}
\end{eqnarray}

By direct computation, we obtain
\begin{eqnarray}\label{2.10}
\Delta_hG&=&\Delta_h(e^{-\lambda(\varphi-\inf_M\varphi)}|\nabla\varphi|_\omega^2)\notag\\
&=&e^{-\lambda(\varphi-\inf_M\varphi)}\Delta_h|\nabla\varphi|_\omega^2+(\Delta_he^{-\lambda(\varphi-\inf_M\varphi)})|\nabla\varphi|_\omega^2\notag\\
&&+2\nabla e^{-\lambda(\varphi-\inf_M\varphi)}\cdot_h\nabla|\nabla\varphi|_\omega^2\notag\\
&=&e^{-\lambda(\varphi-\inf_M\varphi)}[-2ne^F\nabla F\cdot_\omega \nabla\varphi-2tr_h\nabla\chi\cdot_\omega\nabla\varphi\notag\\
&&+2h^{i\bar{j}}g^{k\bar{q}}g^{p\bar{l}}R(\omega)_{i\bar{j}k\bar{l}}\varphi_p\varphi_{\bar{q}}+h^{i\bar{j}}g^{k\bar{l}}(\varphi_{ki}\varphi_{\bar{j}\bar{l}}+\varphi_{k\bar{j}}\varphi_{i\bar{l}})]\notag\\
&&+G(-\lambda ne^F+\lambda tr_h\chi+\lambda^2|\nabla\varphi|_h^2)\notag\\
&&-2\lambda e^{-\lambda(\varphi-\inf_M\varphi)}\nabla\varphi\cdot_h\nabla|\nabla\varphi|_\omega^2.
\end{eqnarray}

Next, we deal with the right-hand-side of \eqref{2.10} in a holomorphic normal coordinate system under K\"ahler metric $\omega$ so that
\[g_{i\bar{j}}=\delta_{ij},\ \chi_{i\bar{j}}=\tau_i\delta_{ij},\ \text{and}\ (\chi_\varphi)_{i\bar{j}}=\rho_i\delta_{ij}\]
for some $\tau_1$,...,$\tau_n>0$ and $\rho_1$,...,$\rho_n>0$. Hence, $h^{i\bar{j}}$ can be diagonalized that
\begin{equation}\label{2.x}
h^{i\bar{i}}=(\chi_\varphi^{i\bar{i}})^2=\frac{1}{\rho_i^2}\leq(\sum_{i=1}^n\frac{1}{\rho_i})^2=(\chi^{i\bar{j}}_\varphi g_{i\bar{j}})^2=n^2e^{2F},
\end{equation}
where we used \eqref{2.5x}.

It follows from Cauchy's inequality and \eqref{2.x} that
\begin{eqnarray}\label{2.20}
&&e^{-\lambda(\varphi-\inf_M\varphi)}(-2ne^F\nabla F\cdot_\omega \nabla\varphi-2tr_h\nabla\chi\cdot_\omega\nabla\varphi)\notag\\
&\geq&(-2ne^F|\nabla F|_\omega-2n^2e^{2F}|\nabla\chi|_\omega)\cdot e^{-\frac{\lambda}{2}(\varphi-\inf_M\varphi)}|\nabla\varphi|_\omega\notag\\
&\geq&-CG^\frac{1}{2}
\end{eqnarray}
for a uniform constant depends only on $M$, $\omega$, $\chi$ and $F$.
\begin{eqnarray*}
-2\lambda\nabla\varphi\cdot_h\nabla|\nabla\varphi|_\omega^2&=&-2\lambda Re(h^{i\bar{j}}\varphi_i(\varphi_k\varphi_{\bar{k}})_{\bar{j}})\notag\\
&=&-2\lambda Re(h^{i\bar{j}}\varphi_i\varphi_{k\bar{j}}\varphi_{\bar{k}}+h^{i\bar{j}}\varphi_i\varphi_k\varphi_{\bar{k}\bar{j}})\notag\\
&\geq&-2\lambda h^{i\bar{j}}\varphi_i\varphi_{\bar{k}}(\chi_\varphi)_{k\bar{j}}+2\lambda h^{i\bar{j}}\varphi_i\varphi_{\bar{k}}\chi_{k\bar{j}}\notag\\
&&-h^{i\bar{j}}\varphi_{ki}\varphi_{\bar{k}\bar{j}}-\lambda^2|\nabla\varphi|_h^2|\nabla\varphi|_\omega^2,
\end{eqnarray*}
where we used Cauchy's inequality.

Note that
\begin{eqnarray*}
-h^{i\bar{j}}\varphi_i\varphi_{\bar{k}}(\chi_\varphi)_{k\bar{j}}&=&-\chi_\varphi^{i\bar{p}}\chi_\varphi^{p\bar{j}}(\chi_\varphi)_{k\bar{j}}\varphi_i\varphi_{\bar{k}}\notag\\
&=&-\chi_\varphi^{i\bar{k}}\varphi_i\varphi_{\bar{k}}\notag\\
&\geq&-(tr_{\chi_\varphi}\omega)|\nabla\varphi|_\omega^2\notag\\
&=&-ne^F|\nabla\varphi|_\omega^2,
\end{eqnarray*}
where we used \eqref{2.5x}, and
\[h^{i\bar{j}}\varphi_i\varphi_{\bar{k}}\chi_{k\bar{j}}=\sum_i\frac{\tau_i}{\rho_i^2}\varphi_i\varphi_{\bar{i}}\geq0,\]
we have
\begin{eqnarray}\label{2.11}
&&-2\lambda\nabla\varphi\cdot_h\nabla|\nabla\varphi|_\omega^2\notag\\
&\geq&-2\lambda ne^F|\nabla\varphi|_\omega^2-h^{i\bar{j}}\varphi_{ki}\varphi_{\bar{k}\bar{j}}-\lambda^2|\nabla\varphi|_h^2|\nabla\varphi|_\omega^2.
\end{eqnarray}

Plugging \eqref{2.20} and \eqref{2.11} into \eqref{2.10}, we have
\begin{eqnarray*}
\Delta_hG&\geq&e^{-\lambda(\varphi-\inf_M\varphi)}(-CG^\frac{1}{2}+2h^{i\bar{j}}g^{k\bar{q}}g^{p\bar{l}}R(\omega)_{i\bar{j}k\bar{l}}\varphi_p\varphi_{\bar{q}}\notag\\
&&+h^{i\bar{j}}g^{k\bar{l}}\varphi_{k\bar{j}}\varphi_{i\bar{l}})+G(-3\lambda ne^F+\lambda tr_h\chi)\notag\\
&\geq&-CG^\frac{1}{2}+(\lambda tr_h\chi-2Ktr_h\omega)G-3\lambda ne^FG.
\end{eqnarray*}
\end{proof}

The ABP maximum principle (see e.g. Lemma 9.3 in \cite{elliptic} of Lemma 5.2 in \cite{chencheng1}) is essential to derive a gradient estimate for \eqref{1.1}. For reader's convenience, we present this result below.
\begin{lemma}[ABP maximum principle]
Let $\Omega\subset\mathbb{R}^d$ be a bounded domain. Suppose $u\in C^2(\Omega)\cap C(\overline{\Omega})$. Denote $U=\sup_\Omega u-\sup_{\partial\Omega}u$. Define
\begin{eqnarray*}
\Gamma^-(u,\Omega)&:=&\{x\in \Omega:u(y)\leq u(x)+\nabla u(x)\cdot(y-x),\notag\\
&&\text{for\ any}\ y\in\Omega\ \text{and}\ |\nabla u(x)|\leq\frac{U}{3diam(\Omega)}\}.
\end{eqnarray*}
Then for some dimensional constant $C_d>0$,
\[U\leq C_d\big(\int_{\Gamma^-(u,\Omega)}det(-D^2u)dx\big)^\frac{1}{d}.\]
In particular, suppose $u$ satisfies $a_{ij}\partial_i\partial_ju\geq f$. Here $a_{ij}$ satisfies the ellipticity condition $a_{ij}\xi_i\xi_j\geq0$. Define $D^*= (deta_{ij})^\frac{1}{d}$. Then there exists another dimensional constant $C'_d>0$ so that
\[U\leq C'_ddiam(\Omega)\|\frac{f^-}{D^*}\|_{L^d(\Omega)}.\]
\end{lemma}

Now we are ready to finish the proof of Theorem \ref{thm1.1}.
\\\\\textbf{Proof of Theorem \ref{thm1.1}:} We apply a trick that Chen-Cheng \cite{chencheng1} and Guo-Phong-Tong \cite{gpt21} used. Let $r$ be chosen so that for any $p\in M$, the geodesic ball $B_\omega(p,r):=\{x\in M|d_\omega(x,p)<r\}$ is contained in a single coordinate neighborhood and under this coordinate, we may identify $B_\omega(p,r)$ with a bounded domain in $\mathbb{R}^{2n}$.

Let $0<\theta\leq\frac{1}{2}$ be a given constant to be determined later, and $\eta:M\rightarrow\mathbb{R}_+$ be a cut-off function such that
\[\eta(x)=
\begin{cases}
1& x\in B_\omega(p,r/2),\\
1-\theta& x\in B_\omega(p,r)\backslash B_\omega(p,3r/4),
\end{cases}
\]
and $\frac{1}{2}\leq1-\theta\leq\eta\leq1$ in the annulus between. Moreover,
\begin{equation}\label{2.12}
|\nabla\eta|_\omega^2\leq\frac{c_0\theta^2}{r^2}\ \text{and}\ |\nabla^2\eta|_\omega\leq\frac{c_0\theta}{r^2}
\end{equation}
for some fixed constant $c_0$ that depends only on $\omega$.

By direct computation, we have
\begin{eqnarray}\label{2.13}
\Delta_h(\eta G^2)&=&\eta\Delta_h G^2+4G\nabla\eta\cdot_h\nabla G+G^2\Delta_h\eta\notag\\
&\geq&\eta(2G\Delta_hG+2|\nabla G|_h^2)-|\nabla G|_h^2\notag\\
&&-4G^2|\nabla\eta|_h^2-G^2|\nabla^2\eta|_\omega\cdot tr_h\omega\notag\\
&\geq&-C\eta G^\frac{3}{2}+2(\lambda tr_h\chi-2Ktr_h\omega)\eta G^2\notag\\
&&-6\lambda ne^F\eta G^2+(2\eta-1) |\nabla G|_h^2\notag\\
&&-4G^2|\nabla\eta|_\omega^2\cdot tr_h\omega-G^2|\nabla^2\eta|_\omega\cdot tr_h\omega\notag\\
&\geq&-C\eta G^{\frac{3}{2}}-6\lambda ne^F\eta G^2\notag\\
&&+G^2tr_h[2\eta(\lambda\chi-2K\omega)-\frac{c_0\theta+4c_0\theta^2}{r^2}\omega],
\end{eqnarray}
where we used Cauchy's inequality in the first inequality, Lemma \ref{lem2.1} in the second and \eqref{2.12} in the last.

Choose $\theta=\min\{\frac{1}{2},\frac{r^2}{2c_0},\frac{r}{4\sqrt{c_0}}\}$. There exists a constant $\lambda$ depending only on  $\omega$ and $\chi$ with
\[2\eta(\lambda\chi-2K\omega)-\frac{c_0\theta+4c_0\theta^2}{r^2}\omega\geq0.\]
Therefore, \eqref{2.13} reduces to
\begin{equation}\label{2.14}
\Delta_h(\eta G^2)\geq-C\eta G^\frac{3}{2}-C\eta G^2.
\end{equation}

Before we applying the ABP maximum principle to \eqref{2.14}, we only need to estimate $det(h^{i\bar{j}})^{-1}$.

Under the holomorphic coordinates mentioned in Lemma \ref{lem2.1}, we have $h^{i\bar{i}}=\frac{1}{\rho_i^2}$ for $i=1,2,...,n$. Moreover,
\begin{eqnarray}\label{2.15}
det(h^{i\bar{j}})^{-1}&=&\big(\prod_{i=1}^n\frac{1}{\rho_i^2}\big)^{-1}=\prod_{i=1}^n\rho_i^2\leq\big(\frac{1}{n}(\sum_{i=1}^n\rho_i)^2\big)^n\notag\\
&=&\frac{1}{n^n}(tr_\omega\chi_\varphi)^{2n}.
\end{eqnarray}


Under the condition of \eqref{1.2}, the main theorem of \cite{yli14} shows that
\[tr_\omega\chi_\varphi\leq C\]
for some uniform constant $C$ depending only on $M$, $\omega$, $\chi$ and $F$.

It follows that
\begin{equation}\label{2.16}
det(h^{i\bar{j}})^{-1}\leq C,
\end{equation}
and
\begin{equation}\label{2.18}
\Delta_\omega\varphi=tr_\omega\chi_\varphi-tr_\omega\chi\leq C.
\end{equation}

Integration by parts gives
\begin{equation}\label{2.17}
\int_MG\omega^n=\int_Me^{-\lambda(\varphi-\inf_M\varphi)}|\nabla\varphi|^2_\omega\omega^n=\frac{1}{\lambda}\int_Me^{-\lambda(\varphi-\inf_M\varphi)}\Delta_\omega\varphi\omega^n\leq C,
\end{equation}
where we used \eqref{2.18} and the fact of $e^{-\lambda(\varphi-\inf_M\varphi)}\leq 1$.

For simplicity, we denote $V:=\sup_{B_\omega(p,r)}(\eta G^2)$. Applying the ABP maximum principle to $\eta G^2$ on $B_\omega(p,r)$, we obtain that
\begin{eqnarray}\label{2.19}
&&V\notag\\
&\leq&\sup_{\partial B_\omega(p,r)}(\eta G^2)+C\cdot r\big(\int_{B_\omega(p,r)}[\eta G^\frac{3}{2}+\eta G^2]^{2n}\omega^n\big)^{\frac{1}{2n}}\notag\\
&\leq&(1-\theta)V+C\cdot \big[ V^\frac{3}{4}+V^{1-\frac{1}{2n}}(\int_{B_\omega(p,r)}G\omega^n)^\frac{1}{2n}\big]\notag\\
&\leq&(1-\theta)V+C\cdot (V^\frac{3}{4}+V^{1-\frac{1}{2n}}),
\end{eqnarray}
where we used \eqref{2.16} in the first inequality and \eqref{2.17} in the last.

We can conclude from \eqref{2.19} that
\[G^2(p)\leq V\leq C.\]
The arbitrary of $p$ implies
\[|\nabla\varphi|_\omega^2\leq C e^{\lambda(\varphi-\inf_M\varphi)}\]
on $M$.$\hfill\Box$
\section{Gradient estimate for parabolic Donaldson's equation} 
\label{sec:3}
We give a proof of Theorem \ref{thm1.2} in this section. Using the notations in Section \ref{sec:2}, we have a similar inequality as in Lemma \ref{lem2.1} for $$H\equiv H(x,t):=\exp\{-\lambda(\varphi(x,t)-\inf_{M\times[0,T]}\varphi)\}|\nabla\varphi|_\omega^2(x,t),$$
where $\lambda$ is a positive constant (depending only on $\omega$ and $\chi$) to be determined later.
\begin{lemma}\label{lem3.1}
There exists a uniform constant $C$ depending only on $M$, $\omega$, $\chi$ and $F$ so that
\begin{eqnarray}\label{3.1}
&&(\Delta_h-n\exp\{F-\frac{\partial\varphi}{\partial t}\}\frac{\partial}{\partial t})H\notag\\
&\geq&-(1+\exp\{F-\frac{\partial\varphi}{\partial t}\})C H^\frac{1}{2}+(\lambda tr_h\chi-2Ktr_h\omega)H\notag\\
&&-3\lambda n\exp\{F-\frac{\partial\varphi}{\partial t}\}H
\end{eqnarray}
on $M\times(0,T]$.
\end{lemma}

\begin{proof}

Calculating in a holomorphic normal coordinate system for the K\"ahler metric $\omega$ at an arbitrary point $x$ as \eqref{2.2} to \eqref{2.5}, we can obtain that
\begin{eqnarray}\label{3.2}
&&(\Delta_h-n\exp\{F-\frac{\partial\varphi}{\partial t}\}\frac{\partial}{\partial t})|\nabla\varphi|_\omega^2\notag\\
&=&2h^{i\bar{j}}R(\omega)_{i\bar{j}k\bar{l}}\varphi_{\bar{k}}\varphi_l+h^{i\bar{j}}(\varphi_{k\bar{j}}\varphi_{\bar{k}i}+\varphi_{ki}\varphi_{\bar{k}\bar{j}})\notag\\
&&+h^{i\bar{j}}(\varphi_{i\bar{j}k}\varphi_{\bar{k}}+\varphi_k\varphi_{i\bar{j}\bar{k}})\notag\\
&&-2n\exp\{F-\frac{\partial\varphi}{\partial t}\}\nabla\frac{\partial\varphi}{\partial t}\cdot_\omega\nabla\varphi.
\end{eqnarray}

It follows from \eqref{1.10} that
\begin{equation}\label{3.3}
n\exp\{F-\frac{\partial\varphi}{\partial t}\}=n\frac{\omega\wedge\chi^{n-1}_\varphi}{\chi^n_\varphi}=tr_{\chi_\varphi}\omega.
\end{equation}

Hence,
\[n\exp\{F-\frac{\partial\varphi}{\partial t}\}\cdot (F_k-\frac{\partial\varphi_k}{\partial t})=-h^{i\bar{j}}\varphi_{i\bar{j}k}-h^{i\bar{j}}\chi_{i\bar{j},k}\]
and
\[n\exp\{F-\frac{\partial\varphi}{\partial t}\}\cdot (F_{\bar{k}}-\frac{\partial\varphi_{\bar{k}}}{\partial t})=-h^{i\bar{j}}\varphi_{i\bar{j}{\bar{k}}}-h^{i\bar{j}}\chi_{i\bar{j},\bar{k}}.\]
Moreover, we can obtain
\begin{eqnarray}\label{3.4}
&&h^{i\bar{j}}(\varphi_{i\bar{j}k}\varphi_{\bar{k}}+\varphi_k\varphi_{i\bar{j}\bar{k}})-2n\exp\{F-\frac{\partial\varphi}{\partial t}\}\nabla\frac{\partial\varphi}{\partial t}\cdot_\omega\nabla\varphi\notag\\
&=&-2n\exp\{F-\frac{\partial\varphi}{\partial t}\}\nabla F\cdot_\omega\nabla\varphi-2tr_h\nabla\chi\cdot_\omega\nabla\varphi.
\end{eqnarray}

Applying \eqref{3.4} to \eqref{3.2}, we have
\begin{eqnarray}\label{3.5}
&&(\Delta_h-n\exp\{F-\frac{\partial\varphi}{\partial t}\}\frac{\partial}{\partial t})|\nabla\varphi|_\omega^2\notag\\
&=&2h^{i\bar{j}}R(\omega)_{i\bar{j}k\bar{l}}\varphi_{\bar{k}}\varphi_l+h^{i\bar{j}}(\varphi_{k\bar{j}}\varphi_{\bar{k}i}+\varphi_{ki}\varphi_{\bar{k}\bar{j}})\notag\\
&&-2n\exp\{F-\frac{\partial\varphi}{\partial t}\}\nabla F\cdot_\omega\nabla\varphi-2tr_h\nabla\chi\cdot_\omega\nabla\varphi.
\end{eqnarray}

The rest computation is follows almost line by line from the proof of Lemma \ref{lem2.1} with minor modifications. Therefore, we can conclude \eqref{3.1}.
\end{proof}

We present the following parabolic version of ABP maximum principle (see \cite{lieb}, Theorem 7.1 or \cite{cc21}, Lemma 3.8) that serves as the main analytic tool in the proof of Theorem \ref{thm1.2}.
\begin{lemma}[Parabolic ABP maximum principle]
Let $\Omega\subset\mathbb{R}^d$ be a bounded domain and $u$ be a smooth function on $\Omega\times[0,T]$, such that
\[b^{ij}\partial_{i}\partial_{j}u-\frac{\partial u}{\partial t}\geq f.\]
Here $b^{ij}$ satisfies the ellipticity condition $b^{ij}\xi_i\xi_j\geq0$. Define $D= (detb^{ij})^\frac{1}{d+1}$. Then there exists a dimensional constant $C_d>0$ so that
\[\sup_{\Omega\times[0,T]}u\leq \sup_{\partial_P(\Omega\times[0,T])}u+C_d(diam(\Omega))^\frac{n}{n+1}\|\frac{f^-}{D}\|_{L^{d+1}(E)}.\]
In the above, $\partial_P(\Omega\times[0,T])$ is the parabolic boundary given by $(\Omega\times\{0\})\cup(\partial\Omega\times[0,T])$, $f^-=\max(-f,0)$ and $E$ is the set of $(x,t)$ on which $\partial_tu\geq0$ and $D_x^2u\leq0$.
\end{lemma}

Similar as the arguments as in the elliptic case, we give a proof of Theorem \ref{thm1.2} in the following.
\\\\\textbf{Proof of Theorem \ref{thm1.2}:} We apply a trick that X. Chen and J. Cheng \cite{cc21} used. Let $(q,s)$ be a maximum point of $H$ on $M\times[0,T]$ and identify the geodesic ball $B_\omega(q,k):=\{x\in M|d_\omega(x,q)<k\}$ with an open domain in $\mathbb{R}^{2n}$ under local coordinates. Let $\beta$ $(0<\beta\leq\frac{1}{2})$ be a given constant to be determined later, and $\xi\in[1-\beta,1]$ be a cut-off function on $M$ satisfying
\[\xi(q)=1,\ \xi=1-\beta\ \text{on}\ M-B_\omega(q,k),\]
and
\[\ |\nabla\xi|^2_\omega\leq\frac{c_1\beta^2}{k^2},\ |\nabla^2\xi|_\omega\leq\frac{c_1\beta}{k^2}\]
for some fixed constant $c_1$ that depends only on $\omega$.

By direct computation, we have
\begin{eqnarray}\label{3.6}
&&(\Delta_h-2n\exp\{F-\frac{\partial\varphi}{\partial t}\}\frac{\partial}{\partial t})(\xi H^2)\notag\\
&=&\xi(\Delta_h-2n\exp\{F-\frac{\partial\varphi}{\partial t}\}) H^2+4H\nabla\xi\cdot_h\nabla H+H^2\Delta_h\xi\notag\\
&\geq&\xi(2H(\Delta_h-2n\exp\{F-\frac{\partial\varphi}{\partial t}\})H+2|\nabla H|_h^2)-|\nabla H|_h^2\notag\\
&&-4H^2|\nabla\xi|_h^2-H^2|\nabla^2\xi|_\omega\cdot tr_h\omega\notag\\
&\geq&-(1+\exp\{F-\frac{\partial\varphi}{\partial t}\})C \xi H^\frac{3}{2}+2(\lambda tr_h\chi-2Ktr_h\omega)\xi H^2\notag\\
&&-6\lambda n\exp\{F-\frac{\partial\varphi}{\partial t}\}\xi H^2+(2\xi-1) |\nabla H|_h^2\notag\\
&&-4H^2|\nabla\xi|_\omega^2\cdot tr_h\omega-H^2|\nabla^2\xi|_\omega\cdot tr_h\omega\notag\\
&\geq&-(1+\exp\{F-\frac{\partial\varphi}{\partial t}\})C \xi H^{\frac{3}{2}}-6\lambda n\exp\{F-\frac{\partial\varphi}{\partial t}\}\xi H^2\notag\\
&&+H^2tr_h[2\xi(\lambda\chi-K\omega)-\frac{4c_1\beta^2+c_1\beta}{k^2}\omega],
\end{eqnarray}
where we used Cauchy's inequality in the first and second inequalities and Lemma \ref{lem3.1} in the second inequality.

Let $\beta=\min\{\frac{1}{2},\frac{k^2}{2c_1},\frac{k}{4\sqrt{c_1}}\}$. We can choose a constant $\lambda$ depending only on $\omega$ and $\chi$ with
\[2\xi(\lambda\chi-K\omega)-\frac{4c_1\beta^2+c_1\beta}{k^2}\omega\geq0.\]
Therefore, \eqref{3.6} reduces to
\begin{eqnarray*}
&&(\Delta_h-2n\exp\{F-\frac{\partial\varphi}{\partial t}\}\frac{\partial}{\partial t})(\xi H^2)\notag\\
&\geq&-(1+\exp\{F-\frac{\partial\varphi}{\partial t}\})C\xi H^{\frac{3}{2}}-6\lambda n\exp\{F-\frac{\partial\varphi}{\partial t}\}\xi H^2,
\end{eqnarray*}
i.e.,
\begin{eqnarray*}
&&(\frac{1}{2n}\exp\{F-\frac{\partial\varphi}{\partial t}\}h^{i\bar{j}}\partial_i\partial_{\bar{j}}-\frac{\partial}{\partial t})(\xi H^2)\notag\\
&\geq&-(\exp\{\frac{\partial\varphi}{\partial t}-F\}+1)C\xi H^{\frac{3}{2}}-3\lambda \xi H^2.
\end{eqnarray*}

Eqt. (27) in W. Sun \cite{sun15} shows that
\begin{equation}\label{3.9}
|\frac{\partial \varphi}{\partial t}|\leq\sup_{M\times\{0\}}|\frac{\partial \varphi}{\partial t}|
\end{equation}
on $M\times[0,T]$. Hence,
\begin{eqnarray}\label{3.7}
&&(\frac{1}{2n}\exp\{F-\frac{\partial\varphi}{\partial t}\}h^{i\bar{j}}\partial_i\partial_{\bar{j}}-\frac{\partial}{\partial t})(\xi H^2)\notag\\
&\geq&-(\exp\{\sup_{M\times\{0\}}|\frac{\partial \varphi}{\partial t}|\}+1)C\xi H^{\frac{3}{2}}-3\lambda \xi H^2.
\end{eqnarray}

Consider in a holomorphic normal coordinate system under K\"ahler metric $\omega$ so that
\[g_{i\bar{j}}=\delta_{ij}\ \text{and}\ (\chi_\varphi)_{i\bar{j}}(\cdot,t)=\rho_i\delta_{ij}(\cdot,t)\]
for any fixed $t\in(0,T]$. Similar as \eqref{2.15}, we have
\[det(h^{i\bar{j}})^{-1}\leq\frac{1}{n^n}(tr_\omega\chi_\varphi)^{2n}=\frac{1}{n^n}(tr_\omega\chi+\Delta_\omega\varphi)^{2n}.\]

From Theorem \ref{thm1.1} and Theorem \ref{thm1.2} in W. Sun \cite{sun15}, we know that
\begin{equation}\label{3.8}
\Delta_\omega \varphi\leq \bar{C}
\end{equation}
on $M\times[0,T]$ for some uniform constant $\bar{C}$ depending only on the initial geometric data and $T$. Therefore,
\begin{equation}\label{3.8}
det(h^{i\bar{j}})^{-1}\leq \bar{C}.
\end{equation}
For any fixed $t\in(0,T]$, integration by parts gives
\begin{eqnarray}\label{3.11}
\int_{M\times[0,T]}G\omega^ndt&=&\int_0^T\int_Me^{-\lambda(\varphi(x,t)-\inf_{M\times[0,T]}\varphi)}|\nabla\varphi|^2_\omega(x,t)\omega^ndt\notag\\
&=&\int_0^T\big(\frac{1}{\lambda}\int_Me^{-\lambda(\varphi(x,t)-\inf_{M\times[0,T]}\varphi)}\Delta_\omega\varphi(x,t)\omega^n\big)dt\notag\\
&\leq& \bar{C}.
\end{eqnarray}



Then we finish the proof by applying the parabolic ABP maximum principle to $\xi H^2$ on $B_\omega(q,k)\times[0,T]$. Note that $H^2(q,s)=\sup_{B_\omega(q,k)\times[0,T]}(\xi H^2)$, we have
\begin{eqnarray}\label{3.12}
H^2(q,s)&\leq&\sup_{\partial_P(B_\omega(q,k)\times[0,T])}(\xi H^2)\notag\\
&&+\bar{C}\cdot k^\frac{2n}{2n+1}\big(\int_{B_\omega(q,k)\times[0,T]}(\xi H^{\frac{3}{2}}+ \xi H^2)^{2n+1}\omega^ndt\big)^{\frac{1}{2n+1}}\notag\\
&\leq&(1-\beta)H^2(q,k)+\bar{C}\cdot \big[ H^\frac{3}{2}(q,s)\notag\\
&&+H^{2-\frac{1}{2n+1}}(q,s)(\int_{B_\omega(q,k)\times[0,T]}H\omega^ndt)^\frac{1}{2n}\big]\notag\\
&\leq&(1-\beta)H^2(q,s)+\bar{C}\cdot (H^\frac{3}{2}(q,s)+H^{2-\frac{1}{2n+1}}(q,s)),
\end{eqnarray}
where we used \eqref{3.8} in the first inequality and \eqref{3.11} in the last.

We can conclude from \eqref{3.12} that
\[(\sup_{M\times[0,T]}H)^2=H^2(q,s)\leq \bar{C}.\]
Moreover,
\[|\nabla\varphi|_\omega^2\leq \bar{C}e^{\lambda(\varphi-\inf_{M\times[0,T]}\varphi)}\]
on $M\times[0,T]$.$\hfill\Box$

\section*{Acknowledgements}
The author would like to thank Prof. Xiaokui Yang and Tao Zheng for helpful discussions.

\end{document}